\documentclass{ajmart}

\issueinfo{2}{3}{September}{2008}

\copyrightinfo{2008}{Aulona Press  \emph{(Albanian J. Math.)}}

\def\issn{{\sc ISSN} 1930-1235: }
\def\issueyear{2008}

\PII{\issn  (\issueyear )}

\pagespan{241}{248}

\usepackage{amsmath, amsbsy,amsfonts,euscript, amssymb,amsbsy,amsthm,latexsym,amsopn,amstext, amsxtra,euscript,amscd}

\usepackage{graphicx}

\usepackage[latin1]{inputenc}

\newtheorem{theorem}{Theorem}

\newtheorem{lemma}{Lemma}
\newtheorem{remark}{Remark}

\theoremstyle{definition}

\newtheorem{exa}{Example}

\def\bP{\mathbb P}
\def\bC{\mathbb C}

\def\Z{\mathbb Z}
\def\X{\mathcal X}

\def\H{\mathcal H}
\def\M{\mathcal M}
\def\L{\mathcal L}

\def\lar{\longrightarrow}

\def\s{\sigma}

\def\X{\mathcal X}
\def\H{\mathcal H}
\def\Hs{\mathcal H_\s}
\def\A{\mathcal A}
\def\D{\Delta}

\def\({\left(}
\def\){\right)}
\def\cO{{\mathcal O}}

\def\l{\lambda}

\def\<{\langle}
\def\>{\rangle}

\def\s{\sigma}

\def\<{\langle}
\def\>{\rangle}
\def\_u{{\mathfrak u}}

\begin{document}

\title{Degree even coverings of elliptic curves by genus 2 curves}

\maketitle

\begin{center} {\sc N. Pjero, M. Ramasaço} \\
\small \vspace*{1ex}
\emph{Dep. of Mathematics,\\
 University of Vlora, Albania\\
 npjero@univlora.edu.al, ramosaco@univlora.edu.al }
\end{center}

\normalsize
\vspace*{2ex}

\begin{center} {\sc T. Shaska} \\
\vspace*{1ex}
\emph{Dep. of Computer Science and Electrical Engineering\\
University of Vlora, Albania}
\end{center}


\begin{abstract}
In this survey we study the genus 2 curves with $(n, n)$-split Jacobian for even $n$.  
\end{abstract}

\section{Introduction}
Let $C$ be a genus 2 curve defined over an algebraically closed field $k$, of characteristic zero. Let $\psi: C
\to E$ be a degree $n$ maximal covering (i.e. does not factor through an isogeny) to an elliptic curve $E$ defined
over $k$. We say that $C$ has a \emph{degree n elliptic subcover}. Degree $n$ elliptic subcovers occur in pairs.
Let $(E, E')$ be such a pair. It is well known that there is an isogeny of degree $n^2$ between the Jacobian $J_C$
of $C$ and the product $E \times E'$. The locus of such $C$, denoted by $\L_n$, is a 2-dimensional algebraic subvariety of the moduli space $\M_2$ of genus two curves and has been the focus of many papers in the last decade; see \cite{Sh1, Sh6, Sh7, SV1, SV2,Fr, FK}.

The space $\L_2$ was studied in Shaska/V\"olklein \cite{SV1}. The space $\L_3$ was studied in \cite{Sh1} were an algebraic description was given as sublocus of $\M_2$. Lately the space $\L_5$ has been studied in detail in \cite{SV2}. The case of even degree has been less studied even though there have been some attempts lately to compute some of the cases for $n=4$; see \cite{deg4}. In this survey we study the genus 2 curves with $(n, n)$-split Jacobian for small $n$. While such curves have been studied by many authors, our approach is simply computational.

\section{Curves of genus 2 with split Jacobians}
Most of the results of this section can be found in \cite{Sh}. Let $C$ and $E$ be curves of genus 2 and 1, respectively. Both are smooth, projective curves defined over $k$, $char(k)=0$. Let $\psi: C \longrightarrow E$ be a covering of degree $n$. From the Riemann-Hurwitz formula, $\sum_{P \in C}\, (e_{\psi}\,(P) -1)=2$ where $e_{\psi}(P)$ is the ramification index of points $P \in C$, under $\psi$. Thus, we have two points of ramification index 2 or one point of ramification index 3. The two points of
ramification index 2 can be in the same fiber or in different fibers. Therefore, we have the following cases of the covering $\psi$:\\

\textbf{Case I:} There are $P_1$, $P_2 \in C$, such that $e_{\psi}({P_1})=e_{\psi}({P_2})=2$, $\psi(P_1) \neq
\psi(P_2)$, and $\forall P \in C\setminus \{P_1,P_2\}$, $e_{\psi}(P)=1$.

\textbf{Case II:} There are $P_1$, $P_2 \in C$, such that $e_{\psi}({P_1})=e_{\psi}({P_2})=2$, $\psi(P_1) =
\psi(P_2)$, and $\forall P \in C\setminus \{P_1,P_2\}$, $e_{\psi}(P)=1$.

\textbf{Case III:} There is $P_1 \in C$ such that $e_{\psi}(P_1)=3$, and $ \forall P \in C \setminus \{P_1\}$,
$e_{\psi}(P)=1$.\\

\noindent In case I (resp. II, III) the cover $\psi$ has 2 (resp. 1) branch points in E.

Denote the hyperelliptic involution of $C$ by $w$. We choose $\mathcal O$ in E such that $w$ restricted to
$E$ is the hyperelliptic involution on $E$. We denote the restriction of $w$ on $E$ by $v$, $v(P)=-P$. Thus,
$\psi \circ w=v \circ \psi$. E[2] denotes the group of 2-torsion points of the elliptic curve E, which are
the points fixed by $v$. The proof of the following two lemmas is straightforward and will be omitted.

\begin{lemma} \label{lem_1}
a) If $Q \in E$, then $\forall P \in \psi^{-1}(Q)$, $w(P) \in \psi^{-1}(-Q)$.

b) For all $P\in C$, $e_\psi(P)=e_\psi\,({w(P)})$.
\end{lemma}

Let $W$ be the set of points in C fixed by $w$. Every curve of genus 2 is given, up to isomorphism, by a binary
sextic, so there are 6 points fixed by the hyperelliptic involution $w$, namely the Weierstrass points of $C$. The
following lemma determines the distribution of the Weierstrass points in fibers of 2-torsion points.

\begin{lemma}\label{lem2} The following hold:
\begin{enumerate}
\item $\psi(W)\subset E[2]$
\item If $n$ is an odd number then
\subitem i) $\psi(W)=E[2]$
\subitem ii) If $ Q \in E[2]$ then \#$(\psi^{-1}(Q) \cap W)=1 \mod (2)$
\item If $n$ is an even number then for all $Q\in E[2]$, \#$(\psi^{-1}(Q) \cap W)=0 \mod (2)$
\end{enumerate}
\end{lemma}

Let $\pi_C: C \lar \bP^1$ and $\pi_E:E \lar \bP^1$ be the natural degree 2 projections. The hyperelliptic
involution permutes the points in the fibers of $\pi_C$ and $\pi_E$. The ramified points of $\pi_C$, $\pi_E$
are respectively points in $W$ and $E[2]$ and their ramification index is 2. There is $\phi:\bP^1 \lar \bP^1$
such that the diagram commutes.
\begin{equation}
\begin{matrix}
C & \buildrel{\pi_C}\over\lar & \bP^1\\
\psi \downarrow & & \downarrow \phi \\
E & \buildrel{\pi_E}\over\lar & \bP^1
\end{matrix}
\end{equation}
Next, we will determine the ramification of induced coverings $\phi:\bP^1 \lar \bP^1$. First we fix some
notation. For a given branch point we will denote the ramification of points in its fiber as follows. Any
point $P$ of ramification index $m$ is denoted by $(m)$. If there are $k$ such points then we write $(m)^k$.
We omit writing symbols for unramified points, in other words $(1)^k$ will not be written. Ramification data
between two branch points will be separated by commas. We denote by $\pi_E (E[2])=\{q_1, \dots , q_4\}$ and
$\pi_C(W)=\{w_1, \dots ,w_6\}$.

\subsubsection{The Case When $n$ is Even}
Let us assume now that $deg(\psi)=n$ is an even number. The following theorem classifies the induced
coverings in this case.

\begin{theorem}\label{thm2}
If $n$ is an even number then the generic case for $\psi: C \lar E$ induce the following three cases for
$\phi: \bP^1 \lar \bP^1$:

\begin{description}
\item[I] $ \left ( (2)^\frac {n-2} 2 , (2)^\frac {n-2} 2 , (2)^\frac {n-2} 2 , (2)^ \frac {n}
2 , (2) \right ) $
\item[II] $ \left ( (2)^\frac {n-4} 2 , (2)^\frac {n-2} 2 , (2)^\frac {n} 2 , (2)^ \frac {n} 2 ,
(2) \right ) $
\item[III] $ \left ( (2)^\frac {n-6} 2 , (2)^\frac {n} 2 , (2)^\frac {n} 2 , (2)^ \frac {n} 2 , (2)
\right ) $
\end{description}
Each of the above cases has the following degenerations (two of the branch points collapse to one)

\begin{description}
\item[I]
\begin{enumerate}
\item $\left ( (2)^\frac {n} 2 , (2)^\frac {n-2} 2 ,
 (2)^\frac {n-2} 2 , (2)^ \frac {n} 2 \right )
$ \item $\left ( (2)^\frac {n-2} 2 , (2)^\frac {n-2} 2 , (4) (2)^\frac {n-6} 2 , (2)^ \frac {n} 2 \right ) $
\item $\left ( (2)^\frac {n-2} 2 , (2)^\frac {n-2} 2 ,
 (2)^\frac {n-2} 2 , (4) (2)^ \frac {n-4} 2 \right )
$ \item $\left ( (3) (2)^\frac {n-4} 2 , (2)^\frac {n-2} 2 ,
 (2)^\frac {n-2} 2 , (2)^ \frac {n} 2 \right )
$
\end{enumerate}
\item[II]
\begin{enumerate}
\item $\left ( (2)^\frac {n-2} 2 , (2)^\frac {n-2} 2 , (2)^\frac {n} 2 , (2)^ \frac {n} 2
\right ) $
\item $\left ( (2)^\frac {n-4} 2 , (2)^\frac {n} 2 , (2)^\frac {n} 2 , (2)^ \frac {n} 2 \right)$
\item $\left ((4) (2)^\frac {n-8} 2, (2)^\frac {n-2} 2 , (2)^\frac {n} 2, (2)^ \frac {n}2\right )$
\item $\left ( (2)^\frac {n-4} 2 , (4) (2)^\frac {n-6} 2 , (2)^\frac {n} 2 , (2)^ \frac {n} 2
\right ) $ \item $\left ( (2)^\frac {n-4} 2 , (2)^\frac {n-2} 2 , (2)^\frac {n-4} 2 , (2)^ \frac {n} 2 \right
) $ \item $\left ((3) (2)^\frac {n-6} 2 , (2)^\frac {n-2} 2 , (4) (2)^\frac {n} 2 , (2)^ \frac {n} 2 \right )
$ \item $\left ( (2)^\frac {n-4} 2 , (3) (2)^\frac {n-4} 2 ,
 (2)^\frac {n} 2 , (2)^ \frac {n} 2 \right )
$
\end{enumerate}
\item[III]
\begin{enumerate}
\item $\left ( (2)^\frac {n-4} 2 , (2)^\frac {n} 2 ,
 (2)^\frac {n} 2 , (4) (2)^ \frac {n} 2 \right )
$ \item $\left ( (2)^\frac {n-6} 2 , (4) (2)^\frac {n-4} 2 ,
 (2)^\frac {n} 2 , (2)^ \frac {n} 2 \right )
$ \item $\left ( (2)^\frac {n} 2 , (2)^\frac {n} 2 ,
 (2)^\frac {n} 2 , (4) (2)^ \frac {n-10} 2 \right )
$ \item $\left ( (3) (2)^\frac {n-8} 2 , (2)^\frac {n} 2 ,
 (2)^\frac {n} 2 , (2)^ \frac {n} 2 \right )
$
\end{enumerate}
\end{description}
\end{theorem}

\proof We skip the details of the proof. \qed
\begin{remark}
The case $n=8$ is the first true generic case when all the subcases occur.
\end{remark}

\subsection{Maximal coverings $\psi:C \lar E$.}
Let $\psi_1:C \lar E_1$ be a covering of degree $n$ from a curve of genus 2 to an elliptic curve. The
covering $\psi_1:C \lar E_1$ is called a \textbf{maximal covering} if it does not factor through a nontrivial
isogeny. A map of algebraic curves $f: X \to Y$ induces maps between their Jacobians $f^*: J_Y \to J_X$ and
$f_*: J_X \to J_Y$. When $f$ is maximal then $f^*$ is injective and $ker (f_*)$ is connected, see \cite{Sh7}
for details.

Let $\psi_1:C \lar E_1$ be a covering as above which is maximal. Then ${\psi^*}_1: E_1 \to J_C$ is injective
and the kernel of $\psi_{1,*}: J_C \to E_1$ is an elliptic curve which we denote by $E_2$; see \cite{FK}. For a fixed Weierstrass point $P \in C$, we can embed $C$ to its Jacobian via
\begin{equation}
\begin{split}
i_P: C & \lar J_C \\
 x & \to [(x)-(P)]
\end{split}
\end{equation}
Let $g: E_2 \to J_C$ be the natural embedding of $E_2$ in $J_C$, then there exists $g_*: J_C \to E_2$. Define
$\psi_2=g_*\circ i_P: C \to E_2$. So we have the following exact sequence
$$ 0 \to E_2 \buildrel{g}\over\lar J_C \buildrel{\psi_{1,*}}\over\lar E_1 \to 0 $$
The dual sequence is also exact
$$ 0 \to E_1 \buildrel{\psi_1^*}\over\lar J_C \buildrel{g_*}\over\lar E_2 \to 0 $$
If $deg (\psi_1)$ is an odd number then the maximal covering $\psi_2: C \to E_2$ is unique up to isomorphism
of elliptic curves. If the cover $\psi_1:C \lar E_1$ is given, and therefore $\phi_1$,
we want to determine $\psi_2:C \lar E_2$ and $\phi_2$. The study of the relation between the ramification
structures of $\phi_1$ and $\phi_2$ provides information in this direction. The following lemma (see
\cite[pg. 160]{FK}) answers this question for the set of Weierstrass points $W=\{P_1, \dots , P_6\}$ of C
when the degree of the cover is odd.

\begin{lemma} Let $\psi_1:C \lar E_1$, be maximal of degree $n$.
Then, the map $\psi_2: C\to E_2$ is a maximal covering of degree $n$. Moreover,
\begin{enumerate}
\item [i) ] if $n$ is odd and ${\cO}_i\in E_i[2]$, $i=1, 2$ are the places such that $\#
(\psi_i^{-1}({\cO }_i)\cap W) = 3$, then $\psi_1^{-1}({\cO }_1)\cap W$ and $\psi_2^{-1}({\cO }_2)\cap W$ form
a disjoint union of $W$.
\item [ii)] if $n$ is even and $Q\in E[2]$, then $\# \left( \psi^{-1}(Q)\right) \cap W= 0$ or 2.
\end{enumerate}
\end{lemma}
The above lemma says that if $\psi$ is maximal of even degree then the corresponding induced covering can
have only type \textbf{I} ramification, see Theorem~\ref{thm2}.
\begin{exa}
Let $\psi: C \to E$ be a degree $n=8$ maximal covering of the elliptic curve $E$ by a genus 2 curve $C$. Then, we have Type I covering as in previous theorem. Hence, the ramification is
\[ \left( (2)^3, (2)^3, (2)^3, (2)^4, (2) \right) \]
This case is the first case which has all its subcases with ramifications as follows:
\begin{equation*}
\begin{split}
\textbf{i)} & \, \left( (2)^4, (2)^3, (2)^3, (2)^4 \right) \\
\textbf{ii)} & \, \left( (2)^3, (2)^3, (4)(2), (2)^4 \right) \\
\textbf{iii)} & \,\left( (2)^3, (2)^3, (2)^3, (4) (2)^2 \right) \\
\textbf{iv)} & \, \left( (3) (2)^2, (2)^3, (2)^3, (2)^4 \right) \\
\end{split}
\end{equation*}
The locus of genus 2 curves in the generic case is a 2-dimensional subvariety of the moduli space $\M_2$. 
It would be interesting to explicitly compute such subvariety since it is the first case which could give some clues to what happens in the general case for even degree. 
\end{exa}

\section{The locus of genus two curves with $(n, n)$ split Jacobians}
In this section we will discuss the Hurwitz spaces of coverings with ramification as in the previous section and
the Humbert spaces of discriminant $n^2$.

\subsection{Hurwitz spaces of covers $\phi : \bP^1 \to \bP^1$}
Two covers $f:X\to \bP^1$ and $f':X'\to \bP^1$ are called \textbf{weakly equivalent} if there is a
homeomorphism $h:X\to X'$ and an analytic automorphism $g$ of $\bP^1$ (i.e., a Moebius transformation) such
that $g\circ f=f'\circ h$. The covers $f$ and $f^\prime$ are called \textbf{equivalent } if the above holds
with $g=1$.

Consider a cover $f:X \to \bP^1$ of degree $n$, with branch points $p_1,...,p_r\in \bP^1$. Pick $p\in \bP^1
\setminus\{p_1,...,p_r\}$, and choose loops $\gamma_i$ around $p_i$ such that $\gamma_1,...,\gamma_r$ is a
standard generating system of the fundamental group $\Gamma:=\pi_1( \bP^1 \setminus\{p_1,...,p_r\},p)$, in
particular, we have $\gamma_1 \cdots \gamma_r=1$. Such a system $\gamma_1,...,\gamma_r$ is called a homotopy basis
of $\bP^1 \setminus\{p_1,...,p_r\}$. The group $\Gamma$ acts on the fiber $f^{-1}(p)$ by path lifting, inducing a
transitive subgroup $G$ of the symmetric group $S_n$ (determined by $f$ up to conjugacy in $S_n$). It is called
the \textbf{monodromy group} of $f$. The images of $\gamma_1,...,\gamma_r$ in $S_n$ form a tuple of permutations
$\s=(\s_1,...,\s_r)$ called a tuple of \textbf{branch cycles} of $f$.

We say a cover $f:X\to\bP^1$ of degree $n$ is of type $\s$ if it has $\s$ as tuple of branch cycles relative to
some homotopy basis of $\bP^1$ minus the branch points of $f$. Let $\Hs$ be the set of weak equivalence classes of
covers of type $\s$. The \textbf{Hurwitz space} $\Hs$ carries a natural structure of an quasiprojective variety.

We have $\Hs=\H_\tau$ if and only if the tuples $\s$, $\tau$ are in the same \textbf{braid orbit} $\mathcal
O_\tau = \mathcal O_\sigma$. In the case of the covers $\phi : \bP^1 \to \bP^1$ from above, the corresponding
braid orbit consists of all tuples in $S_n$ whose cycle type matches the ramification structure of $\phi$.

This and the genus of $\Hs$ in the degenerate cases (see the following table) has been computed in GAP by the
BRAID PACKAGE written by K. Magaard.
\begin{table}[!ht]
\begin{tabular}{c|c|c|c|c|c|c }
deg& Case & cycle type of $\s$ & $\# ( \mathcal O_\s )$ & $G$ & $\dim \Hs$ & genus of $\Hs$ \\
\hline & & && & & \\
8 & & $( 2^3, 2^3, 2^3, 2^4, 2 ) $ & 224 & $S_8$ & 2 & -- \\
 & 1 &  $(2^4, 2^3, 2^3, 2^4 ) $ & 4 & 16 & 1 & 0 \\
 & 2 & $(2^3, 2^3, (4) (2) , 2^4 )$ & 48 & $S_8$ & 1 & 4\\
 & 3 & $(2^3, 2^3, 2^3, (4) (2)^2 ) $ & 96 & $S_8$ & 1 & 16\\
 & 4 & $( (3) 2^2, 2^3, 2^3, 2^4 ) $  & 36 & $S_8$ & 1 & 4 \\
  & & & & & & \\
\end{tabular}
\bigskip
\caption{The length of braid orbits, the order of the group, and the genus of 1-dimensional subspaces for even degree maximal coverings.}
\end{table}

As the reader can imagine even such computations are not easy for higher $n$. It is unclear what are the monodromy groups that appear in all the subcases and the formulas for the lengths of the braid orbits.

\subsection{Humbert surfaces}
Let $\A_2$ denote the moduli space of principally polarized abelian surfaces. It is well known that $\A_2$ is
the quotient of the Siegel upper half space $\mathfrak H_2$ of symmetric complex $2 \times 2$ matrices with
positive definite imaginary part by the action of the symplectic group $Sp_4 (\Z)$.

Let $\D$ be a fixed positive integer and $N_\D$ be the set of matrices
$\tau =
\begin{pmatrix}z_1 & z_2\\
z_2 & z_3
\end{pmatrix}
\in \mathfrak H_2$
such that there exist nonzero integers $a, b, c, d, e $ with the following properties:
\begin{equation}\label{humb}
\begin{split}
& a z_1 + bz_2 + c z_3 + d( z_2^2 - z_1 z_3) + e = 0\\
& \D= b^2 - 4ac - 4de\\
\end{split}
\end{equation}

The \emph{ Humbert surface} $\H_\D$ of discriminant $\D$ is called the image of $N_\D$ under the
canonical map
$$\mathfrak H_2 \to \A_2:= Sp_4( \Z)\setminus{\mathfrak H}_2.$$
It is known that $\H_\D \neq \emptyset$ if and only if $\D > 0$ and $\Delta \equiv 0 \textit { or } 1 \mod 4$. Humbert (1900) studied the zero loci in Eq.~\eqref{humb} and discovered certain relations between points in these spaces and certain plane
configurations of six lines. 

For a genus 2 curve $C$ defined over $\bC$, $[C]$ belongs too $\L_n$ if and only if the isomorphism class $[J_C]
\in \A_2$ of its (principally polarized) Jacobian $J_C$ belongs to the Humbert surface $\H_{n^2}$, viewed as a
subset of the moduli space $\A_2$ of principally polarized abelian surfaces.  There is a one to one correspondence between the points in $\L_n$ and points in $\H_{n^2}$. Thus, we have the map:
\begin{equation}
\begin{split}
& \H_\s \, \, \longrightarrow \, \, \L_n \, \, \longrightarrow \, \, \H_{n^2}\\
([f], (p_1, & \dots , p_r) \to [\X] \to [J_\X]\\
\end{split}
\end{equation}
In particular, every point in $\H_{n^2}$ can be represented by an element of $\mathfrak H_2$ of the form
$$\tau =
\begin{pmatrix}z_1 & \frac 1 n \\
\frac 1 n & z_2
\end{pmatrix}, \qquad z_1, \, z_2 \in \mathfrak H.
$$
There have been many attempts to explicitly describe these Humbert surfaces. For some small discriminant this
has been done by several authors; see \cite{SV1}, \cite{Sh1}. Geometric characterizations of such
spaces for $\D= 4, 8, 9$, and 12 were given by Humbert (1900) in \cite{Hu} and for $\D= 13, 16, 17, 20$, 21
by Birkenhake/Wilhelm (2003).

\section{Computing the locus $\mathcal L_n$ in $\mathcal M_2$}
We take the most general case for maximal coverings of even degree, namely $n$, Type I. The ramification structure of $\phi: \bP^1_x \to \bP^1_z$ is 
\[ \left( (2)^{\frac {n-2} 2}, (2)^{\frac {n-2} 2}, (2)^{\frac {n-2} 2}, (2)^{\frac {n} 2}, (2) \right) \]
We denote the branch points respectively $q_1, \dots , q_5$. Let $q_1=0, q_2=1, q_3=\infty$. The red places in $\bP^1_x$ denote the unramified places and the black places all have ramification index 2. We pick the coordinate $x$ such that it is $x=0, x=1, x=\infty$  in the   unramified places of $\bP^1_z$ and  respectively in the fibers of $0,1,\infty$ as in the picture.

There are exactly $d= \frac {n-2} 2$ places of index 2 in $\phi^{-1} (0)$. Let $P (x)$ denote the polynomial whose roots are exactly these places. Similarly denote by $R(x), Q(x)$ such polynomials for fibers of 1 and $\infty$. The other unramified places in the fibers of 0,1, $\infty$ we denote by $w_4, w_5, w_6$ respectively.

Then, we have
\[ z =  \l  \cdot x \frac {x - w_4} { x- w_6} \cdot  \frac {P^2 (x)} {Q^2 (x)} \]
for some $\l \cdot \in \bC$, $\l \neq 0$. Furthermore, 
\[ z -1 =  \l\cdot (x-1)  \cdot \frac {x - w_5} {x-w_6} \cdot \frac {R^2 (x)} {Q^2 (x)} \]
where $P(x), Q(x), R(x)$ are monic polynomials of degree $d= {\frac {n-2} 2}$ with no multiple roots and no common roots.

\bigskip

\begin{center}
\includegraphics[width=10cm]{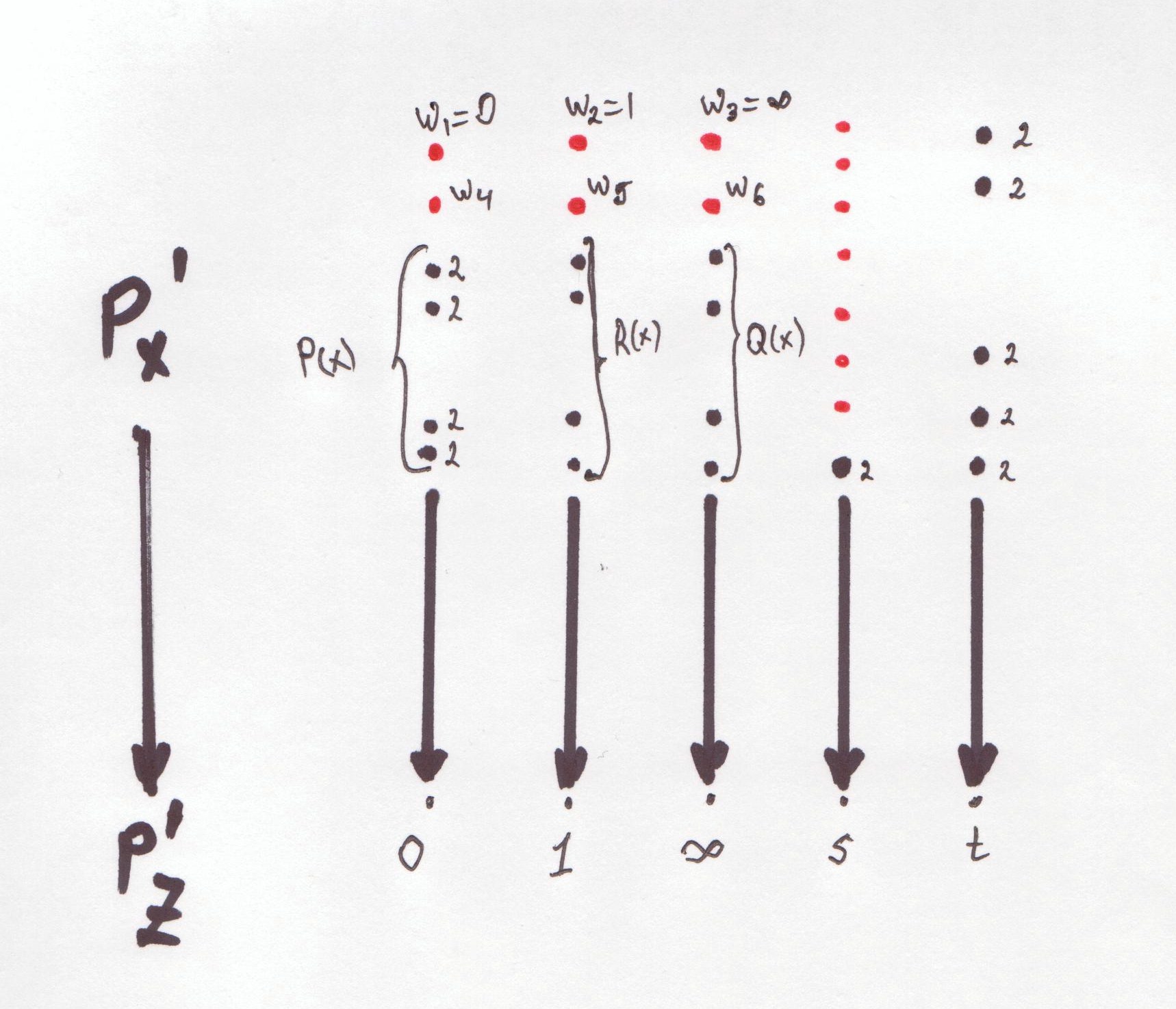}
\end{center}

Substituting for $z$ we get a degree $n$ equation
\[ \l x (x - w_4)   {P^2 (x)}  - ( x- w_6) Q^2 (x)     -  \l \cdot (x-1) (x - w_5) {R^2 (x)} =0 \]
By equaling coefficients of this polynomial with zero we get a nonlinear system of $n+1$ equations. In the same way we get the corresponding equations from the fibers of the other two branch points  $s$ and $t$.  Solving such system would determine also $w_4, w_5, w_6$. The equation of the genus 2 curve $C$ is given by 
\[ y^2= x(x-1)(x-w_4)(x-w_5)(x-w_6)\]

\subsection{Degree 4 covers}
In this section we focus on the case $\deg (\phi) = 4$ (not necessarily maximal). The goal is to determine all ramifications $\s$ and explicitly compute $\L_4 (\s)$.  There is one generic case and one degenerate case in which the ramification of $\deg (\phi) = 4$ applies, as given by the above possible ramification structures.
\begin{enumerate}
\item[i)] $(2,2,2,2^2,2)$ (generic)
\item[ii)] $(2,2,2,4)$ (degenerate)
\end{enumerate}

\subsection{Degenerate Case}
In this case one of the Weierstrass points has ramification index 3, so the cover is totally ramified at this point.

Let the branch points be 0, 1, $\lambda$, and $\infty$, where $\infty$ corresponds to the element of index 4. Then,
above the fibers of 0, 1, $\lambda$ lie two Weierstrass points. The two Weierstrass points above $0$ can be written as
the roots of a quadratic polynomial $x^2 + ax + b$; above $1$, they are the roots of $x^2 +px + q$; and above
$\lambda$, they are the roots of $x^2 + sx + t$. This gives us an equation for the genus 2 curve $C$:
\[ C: y^2 = (x^2 + ax + b)(x^2 +px + q)(x^2 +sx + t). \]
The four branch points of the cover $\phi$ are the 2-torsion points $E[2]$ of the elliptic curve $E$, allowing us to write the elliptic subcover as
\[ E: y^2 = x(x - 1)(x - \lambda). \]
We have the following theorem:

\begin{theorem}
Let $C$ be a genus 2 curve with a degree 4 degenerate elliptic subcover. Then $C$ is isomorphic to the curve given by
\begin{small}
\begin{equation} \label{genus2Degenerate}
\begin{split}
 C: y^2 & = \left(\frac{1-b}{3}+\frac{2}{3} (1-b) x+x^2\right) \left(\frac{1}{12} (b-4) b+\frac{1}{3} (b-4)
 x+x^2\right)\\
 & \left(b-\frac{2}{3} (b+2) x+x^2\right) \\
 E: v^2 & = u (u-1) \left(u - \frac{b^3(4 - b)}{16(b - 1)}\right)
\end{split}
 \end{equation}
\end{small}
where the corresponding discriminants of the right sides must be non-zero. Hence,
\begin{align}
 \D_C: & = b(b-4)(b-2)(b-1) (2+b) \neq 0 \\
 \D_E: & = \frac{(b-4)^2 (b-2)^6 b^6 (b+2)^2}{65536 (b-1)^4} \neq 0.
\end{align} 
and its invariants satisfy 
\begin{tiny}
\begin{equation}\label{loc_deg}
\begin{split}
& 1541086152812576000\,{{ J_2}}^{2}{{ J_4}}^{2}-22835312232360960000 \,{ J_2}\,{J_4}\,{ J_6}+5009676947631\,{{
J_2}}^{6} \\
&-8782271900467200000\,{{ J_6}}^{2} + 1176812184652746480\,{{ J_2}}^{4}{ J_4}+12448207102988800000\,{{
J_4}}^{3}\\
&-3715799948429529600\,{{ J_2}}^{3}{ J_6} =0\\
& 186626560000\,{{J_2}}^{2}{{ J_4}}^{4}+ 138962144767343358744576000000\,{{
J_{10}}}^{2}+\frac {282429536481} {10^4} \,{{J_2}}^{10}\\
&+619923800736 \, {{ J_2}}^{6}{{ J_4}}^{2}-25600000000 \,{{ J_4}}^{5}- \frac {28249152375924} {100} \,{{J_2}}^{8}{J_4}\\
& +266576269949878792320\,{{J_2}}^{5}{J_{10}}-510202022400\, {{J_2}}^{4}{{
J_4}}^{3}\\
&+693067624145203200000\,{ J_2}\,{{ J_4}}^{2}{ J_{10}}+1763516708182388736000\,{{J_2}}^{3}{J_4}\,{J_{10}} =0.
\end{split}
\end{equation}
\end{tiny}
\end{theorem}

\proof  See \cite{deg4}. \qed

\end{document}